\documentclass[11pt]{article} 
\usepackage[letterpaper,hmargin=1in,vmargin=1in]{geometry}
\usepackage{amssymb, amsmath, amsthm, hyperref}

\def\RR{{\mathbb R}}
\def\ZZ{{\mathbb Z}}
\def\CC{{\mathbb C}}
\def\NN{{\mathbb N}}

\def\I{{\sf I}}
\def\P{{\sf P}}
\def\R{{\sf R}}
\def\EE{{\cal E}}
\def\K{{\sf K}}
\def\E{{\rm E_{\pi}}}
\def\var{{\rm Var_{\pi}}}

\newcommand{\edc}{\printindex\end{document}}

\newcommand{\old}[1]{}

\newcommand{\gw}{\omega}
\newcommand{\ga}{\alpha}
\newcommand{\gs}{\sigma}
\newcommand{\ra}[0]{\rightarrow}
\newcommand{\beq}[1]{\begin{equation}\label{#1}}
\newcommand{\enq}[0]{\end{equation}}
\newcommand{\raf}[1]{(\ref{#1})}
\newcommand{\case}[4]{
\left\{ \begin{array}{ll} {#1} & \mbox{#2} \\ {#3}  & \mbox{#4}
\end{array} \right.}

\theoremstyle{plain}
\newtheorem{thm}{Theorem}[section]
\newtheorem{lem}[thm]{Lemma}

\newtheorem{proposition}[thm]{Proposition}

\newtheorem{remark}[thm]{Remark}

\theoremstyle{definition}
\newtheorem{definition}[thm]{Definition}

\begin{document}
\title{Near Optimal Bounds for Collision in Pollard Rho for Discrete Log}

\author{Jeong Han Kim\\
Department of Mathematics\\
Yonsei University\\
Seoul, 120-749 Korea\\
jehkim@yonsei.ac.kr\\
\and
Ravi Montenegro\\
Department of Mathematical Sciences\\
University of Massachusetts Lowell\\
Lowell, MA 01854\\
ravi\_montenegro@uml.edu\\
\and 
Prasad Tetali\\
School of Mathematics and College of Computing\\
Georgia Institute of Technology\\
Atlanta, GA 30332\\
tetali@math.gatech.edu
\thanks{research supported in part by NSF grants DMS 0401239, 0701043}
}

\maketitle

\begin{abstract}
We analyze a fairly standard idealization of Pollard's Rho algorithm for finding the discrete logarithm in  a cyclic group $G$. It is found that, with high probability, a collision occurs 
in  $O(\sqrt{|G|\log |G|\,\log \log |G|})$ steps, not far from the widely conjectured value of $\Theta(\sqrt{|G|})$.  This improves upon a recent result of Miller--Venkatesan which showed an upper bound of $O(\sqrt{|G|}\log^3 |G|)$.
Our proof is based on analyzing an appropriate nonreversible, non-lazy random walk on a discrete cycle of (odd) length $|G|$, and showing that the mixing time of the corresponding walk is $O(\log |G| \log \log |G|)$. 
\end{abstract}

\section{Introduction}
The classical discrete logarithm problem on a cyclic group deals with computing the exponents, given the generator of the group; more precisely, given a generator $x$ of a cyclic group  $G$ and an element $y=x^k$, one would like to compute $k$ efficiently.
Due to its presumed computational difficulty, the problem figures prominently in various cryptosystems, including the Diffie-Hellman key exchange, El Gamal system, and
elliptic curve cryptosystems.  About 30 years ago, J.M. Pollard suggested
algorithms to help solve both factoring large integers \cite{Pol75} and the discrete logarithm problem \cite{Pol78}.
While the algorithms are of much interest in computational number theory and cryptography, there has been very little work on rigorous analyses. We refer the reader to \cite{MV06.1} and other existing literature (e.g., \cite{Tes01,CP05}) for further cryptographic and number-theoretical motivation for the discrete logarithm
problem.

Pollard's Rho algorithm for finding discrete logarithms
is based on a pseudo-random approximation to a Markov chain on a cyclic group $G$.
While there has been no rigorous proof of rapid mixing of the corresponding Markov chain of order $O(\log^c |G|)$ until recently, a proof of mixing of order $O(\log^3|G|)$ steps by a non-trivial argument involving characters and quadratic forms was provided by Miller-Venkatesan \cite{MV06.1}. In addition, they proved that  with high probability the collision time is bounded by $O(\sqrt{|G|} \log^3|G|)$ in the Pollard Rho algorithm. In this paper we improve on this and prove the correct order mixing time bound of a closely related walk, along with a nearly optimal bound on the mixing time of the Pollard Rho algorithm. 

Our first approach will be an elementary proof based on canonical paths which shows the same $O(\log^3|G|)$ mixing time as \cite{MV06.1}. In fact, related methods including log-Sobolev and Spectral profile can show no better than $O(\log^2|G|)$ mixing, still far from the correct bound. As such we then turn to a different method and next show that arguments used to study a related walk by Aldous and Diaconis \cite{AD86.1} and Chung, Diaconis and Graham \cite{CDG87.1} can be modified to apply to this problem. In particular, a strong stationary time is given to show $O(\log |G|\,\,\log \log |G|)$ mixing time when $|G|=2^m-1$ for some $m$, while a Fourier analysis approach can show the same bound for general odd order $|G|$. We then combine this with an improved argument on collision time of a walk, in showing that $O(\sqrt{|G|\log |G|\,\log\log|G|})$ steps suffice until a collision occurs and discrete logarithm is possibly found, not far from the widely conjectured value
of $\Theta(\sqrt{|G|})$.  Finally we observe that our approach is robust enough to allow analysis of other variants of the Pollard Rho algorithm, such as those mentioned in the survey article by Teske \cite{Tes01}; we will include the necessary details in the journal version of this manuscript. A noteworthy remark here is that the walks analyzed in \cite{AD86.1,CDG87.1} and other similar walks studied by Hildebrand \cite{Hil06.1} always double the current position (and then add or subtract one with some probability); the
subtlety in our problem arises from the  original (unaltered) Pollard Rho algorithm demanding that we double only with 1/3 probability. It turns out
that this requires a more careful analysis, since standard comparison-type arguments result in additional $\log p$ factors in the mixing time.

In terms of prior (and not-so-recent) history, we remark  
that Shoup \cite{Shoup97} had shown that any {\em generic}
algorithm which solves (with high probability) the discrete logarithm problem on integers
modulo $p$, must perform at least $\Omega(\sqrt{p})$ group operations, where $p$
is the largest prime dividing $n$. The notion of generic includes among others, Pollard's Rho method and Pohlig-Hellman algorithm (see \cite{Shoup97} for details.) Pollard has shown that if the iterating function $F$ gives perfectly random samples then the expected time until a collision occurs is in fact $O(\sqrt{p})$, but it is not known whether the form of iterating function proposed by Pollard gives a sufficient level of randomness,
and secondly one would like to estimate such a collision time with high probability. 
For one of the  variants of the Pollard Rho algorithm, see
 \cite{HV02.1}, wherein the authors 
replace the squaring step by a walk on a Cayley graph of the group, and obtain
bounds of the form $O(\sqrt{p})$, up to factors of $\log p$.

The paper proceeds as follows. In Section \ref{sec:rho} the Pollard Rho algorithm is introduced, and a relation is shown between collision time and mixing time in separation distance. We then use the canonical path method to bound mixing time of this walk in Section \ref{sec:comparison}. This is followed in Section \ref{sec:SST} by a proof of a near optimal mixing bound in terms of strong stationary times when $|G|=2^m-1$. In Section \ref{sec:fourier}, a Fourier approach is used to show the same bound for the general case. Finally, in conclusion, we 
discuss limitations of various commonly used methods for this problem.


\section{Collision time of Pollard's Rho algorithm} \label{sec:rho}

While the majority of analysis in this paper is devoted to studying the precise mixing time
of a certain nonreversible, non-lazy random walk on  a cycle of odd length,
we first reduce the collision time problem (of  Pollard's Rho discrete logarithm algorithm) to such a mixing time question.  While such a reduction was already described in 
\cite{MV06.1}, our proposition below improves on their idea in yielding a smaller factor (see  below for further clarification). First, let us introduce the algorithm.

Consider a cyclic group $G$ of prime order $p=|G|\neq 2$, and suppose $x$ is a generator, that is $G=\{x^i\}_{i=0}^{p-1}$. Given $y\in G$, the discrete logarithm problem asks us to find $k$ such that $x^k=y$. Pollard suggested an algorithm on $\ZZ_p^{\times}$ based on a random walk and the Birthday Paradox. A common extension of his idea to groups of prime order is to start with a partition of $G$ into sets $S_1$, $S_2$, $S_3$ of roughly equal sizes, and define an iterating function $F:\,G\to G$ by $F(g)=gx$ if $g\in S_1$, $F(g)=g^2$ if $g\in S_2$ and $F(g)=gy=gx^k$ if $g\in S_3$. Then consider the walk $g_{i+1}=F(g_i)$. If this walk passes through the same state twice, say $x^{a+kb}=x^{\alpha+k\beta}$, then $x^{a-\alpha}=x^{k(\beta-b)}$ and so $a-\alpha\equiv k(\beta-b)\mod p$ and $k\equiv (a-\alpha)(\beta-b)^{-1}\mod p$, which determines $k$ unless $\beta\equiv b\mod p$. Hence, if we define a {\em collision} to be the event that the walk passes over the same group element twice, then the first time there is a collision it might be possible to determine the discrete logarithm. 

To estimate the running time until a collision one heuristic is to treat $F$ as if it outputs uniformly random group elements. By the Birthday Paradox if $O(\sqrt{|G|})$ group elements are chosen uniformly at random, then there is a high probability that two of these are the same. However, Teske \cite{Tes98} has given experimental evidence that the time until a collision is slower than what would be expected by a truly random process. We analyze instead the weaker idealization in which it is assumed only that each $g\in G$ is assigned independently and at random to a partition $S_1$, $S_2$ or $S_3$. In this case, although the iterating function $F$ described earlier is deterministic, because the partition of $G$ was randomly chosen then the walk is equivalent to a Markov chain (i.e. a random walk), at least until the walk visits a previously visited state and a collision occurs. The problem is then one of considering a walk on the exponent of $x$, that is a walk $\R$ on the cycle $\ZZ_p$ with transitions $\R(i,i+1)\approx \R(i,i+k)\approx \R(i,2i)\approx 1/3$.

Recall that the event of revisiting an already visited state is called a collision. Our analysis of the time until a collision occurs will be done by examining the rate of convergence of the Markov chain to its stationary distribution $\pi$. The separation distance between a distribution $\sigma$ and stationary distribution $\pi$ is $sep(\sigma,\pi)=\max_{y\in V} 1-\frac{\sigma(y)}{\pi(y)}$. The mixing time of a Markov chain $\P$ with state space $V$ is
$$
\tau_s(\epsilon)=\min\{n:\,\forall x,y\in V,\,\mbox{$1-\frac{\P^n(x,y)}{\pi(y)}$}\leq\epsilon\}\,,
$$
which is the worst-case number of steps required for the separation distance to drop to $\epsilon$. The following result relates $\tau_s(1/2)$ to the time until a collision occurs for any Markov chain $\P$ with uniform distribution on $G$ as the stationary distribution. 

\begin{proposition}
With the above definitions, after
$$
1+\tau_s(1/2) + 2\,\sqrt{2c\,|G|\,\tau_s(1/2)} 
$$
steps, a collision occurs with probability at least $1-e^{-c}$, for any $c>0$.
\end{proposition}

\begin{proof}
Let $S$ denote the first $\left\lceil \sqrt{2c\,|G|\,\tau_s(1/2)}\right\rceil$ states visited by the walk. If two of these states are the same then a collision has occurred, so assume all states are distinct. Even if we only check for collisions every $\tau_s(1/2)$ steps, the chance that no collision occurs in the next $t\tau_s(1/2)$ steps (so consider $t$ semi-random states) is then at most
$$
\left(1-\frac 12\frac{|S|}{|G|}\right)^t
 \leq \left(1-\sqrt{\frac{c\,\tau_s(1/2)}{2|G|}}\right)^t
 \leq e^{-t\,\sqrt{\frac{c\,\tau_s(1/2)}{2|G|}}}\,.
$$
When $t=\left\lceil\sqrt{\frac{2c |G|}{\tau_s(1/2)}}\right\rceil$, this is at most $e^{-c}$, as desired, and so at most
$$
\left\lceil \sqrt{2c\,|G|\,\tau_s(1/2)}\right\rceil
 + \left\lceil\sqrt{\frac{2c |G|}{\tau_s(1/2)}}\right\rceil\,\tau_s(1/2)
$$
steps are required for a collision to occur with probability at least $1-e^{-c}$.
\end{proof}


\begin{remark}
By assuming each $g\in G$ is assigned independently and at random to a partition we have eliminated one of the key features of the Pollard Rho algorithm, space efficiency. However, if the partitions are given by a hash function $f:(G,p)\to\{1,2,3\}$ which is sufficiently pseudo-random then we might expect behavior similar to the model with random partitions.
\end{remark}

\begin{remark}
The analysis can easily be extended to the case when there are many partitions, of varying sizes, each with their own transition rule, as long as a partition occupying at least a constant fraction of the space corresponds to $g\to gx$ and another such partition corresponds to $g\to g^2$.
\end{remark}

Throughout the analysis in the following sections, we assume that the size $p$ of the cycle $\ZZ_p$ (on which the random walk is performed) is odd. Indeed there is a standard reduction -- see \cite{Pom} for a very readable account and also a classical reference \cite{PohHel78} -- justifying the fact that it suffices to study the discrete logarithm problem on cyclic groups of {\em prime} order.


\section{Canonical Paths} \label{sec:comparison}

Perhaps the most widely used approach to bounding mixing times is the method of canonical paths. However, this method has been used primarily for walks which are either lazy or reversible, and usually both. The Pollard Rho walk is neither, but as we will now see, it is still possible to apply the canonical path method. 

Canonical path methods rely on studying the spectral gap:

\begin{definition}
Given Markov chain $\P$ on state space $V$ the spectral gap $\lambda=\lambda_{\P}$ is defined by
$$
\lambda_{\P} = \inf_{\substack{f:V\to\RR,\\ \var(f)\neq 0}} \frac{\EE_{\P}(f,f)}{\var(f)},
$$
with $\var(f) = \E f^2 - (\E f)^2$ and Dirichlet form
$$
\EE_{\P}(f,f) = \frac 12\,\sum_{x,y\in V}(f(x)-f(y))^2\pi(x)\P(x,y)\,.
$$
\end{definition}

Fill \cite{Fill91.1}, building on work of Mihail \cite{Mih89.1}, showed a bound on the mixing time. 

\begin{thm} \label{thm:Fill}
The mixing time of a finite Markov chain $\P$ on state space $V$ 
is at worst
$$
\tau_s(\epsilon) \leq \left\lceil \frac{1}{\lambda_{\P\P^*}}\,\log\frac{1}{\epsilon\pi_0}\right\rceil\,,
$$
where $\pi_0=\min_{x\in V}\pi(x)$ and the time reversal $\P^*$ is given by $\P^*(x,y)=\frac{\pi(y)\P(y,x)}{\pi(x)}$. 
\end{thm}

\noindent
One of the more common ways of bounding the  spectral gap is via canonical paths \cite{Sin92.1}. 

\begin{thm} \label{thm:paths}
Consider a finite Markov chain $\P$ on state space $V$. For every $x,\, y\in V$, $x\neq y$, define a path $\gamma_{xy}$ from $x$ to $y$ along edges of $\P$ (i.e. $\gamma_{xy} \subset E=\{(a,b)\in V\times V:\,\P(a,b)>0\}$). Then
$$
\lambda \geq \Bigl(\max_{(a,b)\in E} \frac{1}{\pi(a)\P(a,b)}\, \sum_{\substack{(x,y):\,x\neq y,\\ (a,b)\in\gamma_{xy}}} \pi(x)\pi(y)|\gamma_{xy}|\Bigr)^{-1}\,.
$$
\end{thm}

It suffices to bound the mixing time of the walk $\R^2$, because the mixing time of $\R$ is at most twice this.

\begin{lem} \label{lem:comparison}
Let $\K(i,2i)=\K(i,2i-1)=1/2$ be a walk on the odd cycle $\ZZ_p$. The Pollard Rho walk $\R$ satisfies
$$
\lambda_{\R^2(\R^2)^*} \geq \frac{2}{81}\,\lambda_{\K}\,.
$$
\end{lem}

\begin{proof}
Observe that, from the definition of spectral gap,
$$
\forall i\neq j:\,\R^2(\R^2)^*(i,j)\geq c \ \K(i,j)
\ \ \Rightarrow\ \ 
\lambda_{\R^2(\R^2)^*} \geq c \ \lambda_{\K} .
$$
Now, $\K(i,j)\neq 0$ only if $j=2i-1$ or $j=2i$, so it suffices to consider these transitions:
\begin{eqnarray*}
\lefteqn{\R^2(\R^2)^*(i,2i-1)} \\
 &\geq& \R(i,2i)\R(2i,2i+1)\R^*(2i+1,2i)\R^*(2i,2i-1) \\
 &\geq& \frac{2}{81}\,\K(i,2i-1)\,, \\
\lefteqn{\R^2(\R^2)^*(i,2i)} \\
 &\geq& \R(i,i+1)\R(i+1,2i+2) \\ 
      &&\R^*(2i+2,2i+1)\R^*(2i+1,2i) \\
 &\geq& \frac{2}{81}\,\K(i,2i)\,.
\end{eqnarray*}
Here we are using the fact that $(\R^2)^*=(\R^*)^2$.
\end{proof}

\begin{lem} \label{lem:spectral}
Let $\K(i,2i)=\K(i,2i-1)=1/2$ be a walk on the odd cycle $\ZZ_p$. Then,
$$
\lambda_{\K} \geq \frac{1}{2 \left(\lceil \log_2 p\rceil\right)^2}\,.
$$
\end{lem}

\begin{proof}
Suppose $x,y\in V$ and let $n=\lceil \log_2 p\rceil$. To construct a path from $x$ to $y$, let $x_0=x$ and consider all possible paths of length $n$, i.e. $x=x_0\to x_1\to\cdots\to x_n$ with $x_i=2x_{i-1}-c_i$ and $c_i\in\{0,1\}$. Then
\begin{equation} \label{eqn:sum}
x_n \equiv 2^n x_0 - \sum_{i=1}^n 2^{n-i}\,c_i \mod p\,.
\end{equation}
Each value in $\{0,1,\ldots,2^n-1\}$ can be written in exactly one way as a sum $\sum_{i=1}^n 2^{n-i}\,c_i$, and so there are either one or two possible paths from $x_0$ to $x_n=y$. Pick one as the canonical path $\gamma_{xy}$.

To apply Theorem \ref{thm:paths}, fix edge $(a,b)$ with $\K(a,b)>0$ and suppose that $(a,b)$ is the $i$-th edge in path $\gamma_{xy}$. Then $x\in\{z:\,\K^{i-1}(z,a)>0\}$ and $y\in\{z:\,\K^{n-i}(b,z)>0\}$, and so there are at most $|\{z:\,\K^{i-1}(z,a)>0\}|\times|\{z:\,\K^{n-i}(b,z)>0\}|\leq 2^{i-1}\times 2^{n-i}=2^{n-1}< p$ such paths. There are $n=\lceil \log_2 p\rceil$ possible values of $i$, so there are at most $p\times\lceil \log_2 p\rceil$ paths through this edge, each of length $|\gamma_{xy}|=\lceil \log_2 p\rceil$. Theorem \ref{thm:paths} completes the proof.
\end{proof}

It follows from Theorem \ref{thm:Fill}, Lemma \ref{lem:comparison} and Lemma \ref{lem:spectral} that $\tau_s(\epsilon)=O((\log p)^2\log(p/\epsilon))$.

The preceding argument was based on the observation that studying $\R^m$ for some $m>1$ may be easier than studying the walk $\R$ directly, as was done in \cite{MV06.1}. In Concluding Remarks we sketch an argument for why this approach cannot be used to show better than $\tau_s(1/2)=O((\log p)^2/m)$ for the Pollard Rho walk.


\section{Separation via Stopping Time} \label{sec:SST}

While canonical paths are much more widely used for bounding mixing time, the most direct route for bounding the separation distance is via a strong stationary time. 

\begin{definition}
A stopping time for a random walk $\{Y_i\}_{i=0}^{\infty}$ is a random variable $T\in{\mathbb N}$ such that the event $\{T\le t\}$ depends only on $Y_0,\,Y_1,\,\ldots,\,Y_t$. A stopping time $T$ is a strong stationary time with stationary distribution $\pi$,  if
$$
\forall y,\,\Pr[Y_t=y|T\le t]=\pi(y)\,.
$$
\end{definition}

The key point here is that 
\begin{eqnarray*}
\Pr[Y_t=y]&\geq& \Pr[T\leq t]\,\Pr[Y_t=y|T\leq t]\\
  &=&\pi(y)\Pr[T\leq t]
\end{eqnarray*}
and so 
$$
sep(\P^t(x,\cdot),\pi)=\max_{y\in V}\ 1-\frac{\P^t(x,y)}{\pi(y)} \leq \Pr[T>t]\,.
$$

We first consider the case $p=2^m-1$.
The construction can be thought of as an extension of the approach of Aldous and Diaconis \cite{AD86.1}. 
In the following section a similar bound for general odd $p$ will be shown.

\begin{thm}
If $p=2^m-1$ then the Pollard Rho walk on the cycle $\ZZ_p$ has mixing time 
$$
\tau_s(1/2) = O(\log p\,\log \log p)\,.
$$
\end{thm}

\begin{proof}
The key to the proof will be to reduce the problem to one of constructing a strong stationary time for a walk with transitions $i\to 2i+0$ and $i\to 2i+1$, each with equal probability.

Let us refer to the three types of moves that the Pollard Rho random walk makes,
namely $(i, i+1), (i,i+k)$, and $(i, 2i)$, as moves of Type 1, Type 2, and Type 3, 
respectively. In general, let the random walk be denoted by $Y_0, Y_1, Y_2, \ldots$\,, with 
$Y_t$ indicating the position of the walk (modulo $p$) at time $t\ge 0$.
\vspace{2ex}

\noindent
{\em Define new random variables $T_i$, $b_i$ and $X_i$:}
Let $X_0=0$ and $T_0=0$. Let $T_1$ be the first
time, after time 0, that the walk makes a move of Type 3.
Let $b_1 = Y_{T_1 - 1} - Y_{T_0}$ (i.e., the ground covered, mod $p$,  only using
consecutive moves of Types 1 and 2.)
Let $X_1 = 2X_0 + b_1$.  (Thus $X_1 = Y_{T_1 - 1}-Y_{T_0}$.)
More generally, let $T_i$ be the first time, since $T_{i-1}$, that a move of Type 3 happens. 
Let $b_i = Y_{T_i - 1} - Y_{T_{i-1}}$,  and
let $X_i = 2X_{i-1} + b_i = Y_{T_i-1}-Y_{T_0}$.
Observe that $T_i$, for each $i$,  is a valid stopping time.
\vspace{2ex}

\noindent
{\em Auxiliary Randomness:}
For the sake of the analysis, we generate the above random walk using an auxiliary
random process: at each time step $t\ge 0$, we generate an integer $R_t$ 
uniformly at random from the set
of integers [1..9]. The integers 1,2,3 are associated with (or interpreted as)  a move of Type 1, and the integers 4,5,6 with  Type 2, and finally the integers 7,8,9 with a move of Type 3.
\vspace{2ex}

\noindent
{\em Define History:}
To keep the independence of random variables transparent, it is best to associate
{\em history} vectors $H_i$,  with the random walk as follows. The entries of the history
vector are from [1..9]; every time a doubling move happens, we stop the current
history vector (after recording the current $R_{T_i}$ value), and start growing a new 
vector. Thus $H_1 = (R_1, ..., R_{T_1})$, and in general,
$H_i = (R_{T_{i-1}+1}, ..., R_{T_i})$, and note that the history vector
always ends in a 7, 8, or a 9, since those are identified with a Type 3 move.
\vspace{2ex}

\noindent
{\em Special history vectors: }
We call certain history vectors of length one or two, 
as special: $H$ is special if $H=(7)$ or $H=(a,b)$, where $a \in \{1,2,3\}$, and $b \in \{7,8,9\}$. 
Note that given the history vectors and $Y_0$, all other (random) variables, $Y_i, b_i, T_i, X_i$, 
can  (uniquely) be determined.  Moreover, if a history $H_i$ were special, then it implies that the corresponding $b_i$ equals 0 or 1, depending on whether $H_i=(7)$ or
$(a,b)$, respectively; in the latter case, $a$ being 1, 2, or 3 implies that a move of Type 1 
took place before the doubling, and hence the ground covered is simply $+1$.  
\vspace{2ex}

\noindent
This completes the set up.
\vspace{2ex}

\noindent
{\em The Actual Analysis:}

Let $s=rm$. (Recall that $m=\log_2(p+1)$; we will choose later $r=c \log_2 m$, for
$c>0$ a suitable constant.) Consider
\[ X_s = 2^{s-1} b_1 + 2^{s-2} b_2 + \cdots + 2^0 b_s\,,\] 
which may be rewritten (using modulo $p$) as
\begin{eqnarray*}
X_s &=& 2^{m-1} (b_1 + b_{m+1} +b_{2m+1} + ... + b_{(r-1)m+1}) \\
       &&      + 2^{m-2} (b_2 + b_{m+2} + b_{2m+2} + ... + b_{(r-1)m+2}) \\
       &&      +  \cdots + 2^0(b_{m} + b_{2m} + ... + b_{rm}).
\end{eqnarray*}

In other words, if we refer to each set of terms inside the parentheses as a Block then
there are $m$ Blocks, each associated with $2^i$ for $i=m-1,m-2, ..., 0$.
\vspace{2ex}

\noindent
{\em Define Auxiliary random variables using Special History:}
Recall that each history vector $H_i$ produces a $b_i$.
Let $C_1 = b_{jm+1}$, where $j \in \{0,1, ..., (r-1)\}$ is the first (smallest) index 
in Block 1 such that $b_{jm+1}$ comes from a special history. 
More generally, for $i=1, 2, ..., m-1$,
let $C_{i} = b_{jm+i}$, where $j \in \{0,1, ..., (r-1)\}$ is the first (smallest) $j$ such that $b_{jm+i}$ comes  from a special history. 
If no such j were present (which is possible, since there need not be any occurrences of special history in the corresponding interval), then denote such a $C_i$ to be 
$\infty$ (or undefined.)

By the remarks above, each $C_i$ (once defined) is 0 or 1, and moreover each $C_i$ is  an independent
(of all the other $C_j$'s) Bernoulli trial, since the corresponding $b_i$'s are mutually independent.
Then we may rewrite $X_s$ as follows:
\begin{eqnarray*}
\lefteqn{X_s = 2^{m-1}C_1 + 2^{m-2}C_2 + \cdots + 2^0 C_m } \\
    &&       + 2^{m-1} ({\rm Rest}_1) + 2^{m-2} ({\rm Rest}_2) + \cdots + ({\rm Rest}_m),
\end{eqnarray*}
where ${\rm Rest}_i$ is the sum (over $j$) of $b_{jm+i}$ minus the 
special $b$ that became $C_i$.
\vspace{2ex}

\noindent
{\em The Basic Dyadic Randomness argument from \cite{AD86.1}:}
What is relevant or important here is that if all $C_i$ are defined then 
\begin{eqnarray*}
X_s &=& 2^{m-1}C_1  + 2^{m-2}C_2  + \cdots + 2^0 C_m  + {\rm REST} \\
&=:& S_m + {\rm REST},
\end{eqnarray*}
where, as we will see shortly, the first part ($S_m$) randomizes $X_s$ so that the {\rm REST} will not matter; more formally, if $S_m\neq\infty$ then
\begin{eqnarray*}
\lefteqn{\Pr[X_s = w]} \\
&=& \sum_R  \Pr[{\rm REST} = R] \Pr[S_m + R = w | {\rm REST}=R],
\end{eqnarray*}
and $\Pr[S_m =  w-R | {\rm REST}=R] = 1/2^m$ except that if $w-R=0$ then $\Pr[S_m=0|{\rm REST}=R]=2/2^m$ (when all $C_i$'s are $0$ or all $1$). This holds even if we condition on $T_s=constant$ as well.

Consider the stopping rule $T$ such that $T=T_{rm}$ for the first $r$ for which all $C_i$ are well-defined, except that if all $C_i$ are $1$ then set all $C_i$ to undefined and begin the process again. The above shows that this is a strong stationary time for the Rho walk. 

It remains to bound the time until all $C_i$ are well-defined, which (by coupon-collector intuition) should be roughly $m \log m$: Observe that for a fixed $i$, 
\begin{equation} \label{eqn:undefined}
\Pr[C_i = \infty]  = (1 - 2/9)^r\,.
\end{equation} 
since $\Pr[appropriate\ history\ H_i\ is\ special] = 1/9 + 1/9$, and each of 
the $r$ (independent!)  possibilities in the $i$th Block ought to have been unsuccessful.

So as long as $r > (1+\delta)(\log m)/\log(9/7)$, the probability in \eqref{eqn:undefined} is at most
$1/m^{(1+\delta)}$.  Hence, $\Pr[all\ C_i \in \{0,1\}] \ge 1 - m/m^{1+\delta} = 1- m^{-\delta}$, and so $\Pr[all\ C_i \in \{0,1\} \cap not\ all\ C_i=1] \ge (1-1/2^m)(1-1/m^{\delta})$. For $s=rm$ with $r=\lceil 3(\log m)/\log(9/7)\rceil$, we have
$$
\Pr[T\leq T_s]\geq(1-1/2^m)(1-1/m^2)\,.
$$

A Type 3 move occurs on average every $3$ steps and so by Markov's inequality since $E[T_s]=3s$ then $\Pr[T_s>9s]<1/3$. Hence, if $k=9m\lceil\frac{3\log m}{\log(9/7)}\rceil$ then
\begin{eqnarray*}
\Pr[T\leq k] &\geq& \Pr[T\leq T_s]-\Pr[T_s>k] \\
  &\geq& (1-1/2^m)(1-1/m^2)-1/3 > 1/2\,.
\end{eqnarray*}
\end{proof}

In Concluding Remarks it will be shown that the mixing bound $\tau_s(1/2)=O(\log p\log\log p)$ found here for $X_s$ is of the correct order.


\section{Fourier Analysis} \label{sec:fourier}

We now turn to the general case of $p$ odd, where we work with Fourier analysis. The construction can be thought of as an application of the ideas of Chung, Diaconis and Graham \cite{CDG87.1}.

In the previous two sections a bound on mixing time of the Rho walk was used to derive a bound on collision time. This time we consider a ``block walk'' in which a single step corresponds to a Rho walk truncated on an $i\to 2i$ step, i.e. a walk $Z_i$ where $Z_i=2(Z_{i-1}+b_i)$ with $b_i$ defined as in Section \ref{sec:SST}.  Note that in $t$ steps of the Rho walk the expected number of block steps is $t/3$, and Chebyshev's Inequality shows that $Prob(\# block\ steps<t/4) \leq \frac{32}{t}$. Hence, if $\tau_s^b(1/2)$ denotes the mixing time of the block walk, then in $t=4(1+\tau_s^b(1/2) + 2\,\sqrt{2c\,|G|\,\tau_s^b(1/2)})$ steps of the Rho walk a collision occurs with probability at least $1-e^{-c}-Prob(\# block\ steps<t/4)\geq 1-e^{-c}-\frac{32}{t}$.

To bound mixing time of the block walk it suffices to show that for large enough $s$ the distribution $\nu_s $ of
\[ X_s =  2^{s-1} b_1 + \cdots + b_{s}\,\]
is close to the uniform distribution $u$. More precisely, we will show
that
\begin{equation} \label{eqn:L2} 
p \sum_{j=0}^{p-1} (\nu_s (j) - u(j))^2 \leq 2 \left( (
1+\xi^{2\lfloor s/m \rfloor})^{m-1} -1 \right),
\end{equation}
where $\nu_s (j) = \Pr[ X_s = j]$, $\xi = 1- \frac{4-\sqrt{10}}{9}$, and $m$
satisfies $2^{m-1} < p < 2^m$. This suffices to bound the separation distance, as shown in Remark \ref{rmk:L2-to-sep} at the end of the section.

The proof uses the standard Fourier transform and the Plancherel identity:
For any complex-valued function $f$ on $\ZZ_p$ and $\gw= e^{2\pi i/p}$, 
recall that  the Fourier transform $\hat{f}:\ZZ_p\ra  C$ is given by
$\displaystyle \hat{f} (\ell) = \sum_{j=0}^{p-1} \gw^{\ell j } f(j)$, and the
Plancherel identity asserts that
$$
p \sum_{j=0}^{p-1} |f(j)|^2= \sum_{j=0}^{p-1} |\hat{f}(j)|^2\,.
$$

\noindent
For the distribution $\mu$ of a $\ZZ_p$-valued random variable $X$, its
Fourier transform is
$$
\hat{\mu} (\ell) 
= \sum_{j=0}^{p-1} \gw^{\ell j } \mu (j)
=  E[ \gw^{\ell X}]. 
$$ 
Thus, for the distributions $\mu_{_1}, \mu_{_2}$ of two independent random variables $Y_1, Y_2$, the distribution $\nu$ of
$X:=Y_1+Y_2$ has the Fourier transform $\hat{\nu}=\hat{\mu}_{_1}\hat{\mu}_{_2}$, since

\begin{eqnarray*}
\hat{\nu} (\ell) 
  &=& E[ \gw^{\ell X}] 
   =  E[ \gw^{\ell (Y_1+Y_2)}] \\
  &=& E[ \gw^{\ell Y_1}]E[ \gw^{\ell Y_2}] 
   =  \hat{\mu}_{_1} (\ell)\hat{\mu}_{_2}(\ell).
\end{eqnarray*}
Generally, the distribution $\nu$ of $X:=Y_1 +\cdots + Y_s$
with independent $Y_i$'s has the Fourier transform
$\hat{\nu}=\prod_{j=1}^s  \hat{\mu}_{_j}$.
Moreover, for the uniform distribution $u$, it is easy to check that
$$ 
\hat{u} (\ell) 
= \case{1}{if $\ell=0$,}
       {0}{otherwise.}
$$
As the random variables $2^{j}b_{s-j}$'s are independent,
$\hat{\nu}_s = \prod_{j=0}^{s-1} \hat{\mu_j} $, where $\mu_j$ are
the distributions of $2^{j}b_{s-j}$. The linearity of the Fourier transform and
$\hat{\nu}_s (0) = E[1]=1$ yield
$$ 
\widehat{\nu_s - u} (\ell) 
  =  \hat{\nu}_s (\ell)- \hat{u}(\ell)
  =  \case{0}{if $\ell=0$}{\prod_{j=0}^{s-1} \hat{\mu_j}(\ell)}{otherwise.}
$$
By Plancherel's identity,  it is enough to show that

\begin{lem}
$$
\sum_{\ell=1}^{p-1} \Big|\prod_{j=0}^{s-1} \hat{\mu_j}
(\ell) \Big|^2 \leq 2 \left( (1+\xi^{2\lfloor s/m \rfloor })^{m-1}
-1 \right). 
$$
\end{lem}

\begin{proof}
Let $A_j$ be the event that $b_{s-j} = 0 ~{\rm or} ~1$. Then,
\begin{eqnarray*}
\hat{\mu_j}(\ell) 
  &=& E [  \gw^{\ell 2^{j} b_{s-j}}] \\
  &=& \Pr[ b_{s-j}=0] + \Pr[b_{s-j}=1] \gw^{\ell 2^{j} } \\
  && + \Pr[ \bar{A}_j ]E [  \gw^{\ell 2^{j} b_{s-j}}| \bar{A}_j], 
\end{eqnarray*}
and, for $x:=\Pr[b_{s-j}=0]$ and $y:= \Pr[b_{s-j}=1]$,
\begin{eqnarray*}
|\hat{\mu_j}(\ell)| 
  &\leq& | x + y \gw^{\ell 2^{j} }| + (1-x-y)| E [  \gw^{\ell 2^{j} b_{s-j}}| \bar{A}_j]| \\
  &\leq& | x + y \gw^{\ell 2^{j} }| + 1-x-y .
\end{eqnarray*}
Notice that
\begin{eqnarray*}   
| x + y \gw^{\ell 2^{j} }|^2 
  &=& ( x+ y \cos \mbox{$\frac{2\pi\ell 2^j}{p}$} )^2+ y^2 \sin^2 \mbox{$\frac{2\pi\ell 2^j}{p}$} \\ 
  &=&   x^2 + y^2 + 2 xy \cos \mbox{$\frac{2\pi\ell 2^j}{p}$}. 
\end{eqnarray*}
If $\cos \mbox{$\frac{2\pi\ell 2^j}{p}$}\leq 0$, then
\begin{eqnarray*}
|\hat{\mu}_j (\ell)| 
  &\leq& (x^2 +y^2)^{1/2}+1-x-y \\
  &=& 1-(x+y-(x^2 +y^2)^{1/2})
\end{eqnarray*}
Since $x= \Pr[b_ {s-j}=0] \geq 1/3$ and $y= \Pr[b_{s-j}=1]\geq 1/9$, 
it is easy to see that $x+y-(x^2 +y^2)^{1/2}$ has its minimum when $x=1/3$ and $y=1/9$.
(For both partial derivatives are positive.) Hence,
$$
|\hat{ \mu}_j (\ell)| \leq \xi=1- \frac{4-\sqrt{10}}{9}, 
  ~~~{\rm provided} ~  \cos \mbox{$\frac{2\pi\ell 2^j}{p}$}
  \leq 0.  
$$
If $ \cos \mbox{$\frac{2\pi\ell 2^j}{p}$}> 0$, we use the trivial bound 
$\hat{\mu}_j (\ell) = E [  \gw^{\ell 2^{j} b_{s-j}}] \leq 1$.

For $\ell=1,..., p-1$, let $\phi_s  (\ell) $ be the number of $j=0,..., s-1$ such that 
$\cos \mbox{$\frac{2\pi\ell 2^j}{p}$}\leq 0$. Then
  \beq{phibd} \prod_{j=0}^{s-1} | \hat{\mu}_j (\ell)|  \leq \xi^{\phi_s (\ell)}.
  \enq
To estimate $\phi_s (\ell)$, we consider the binary expansion of
$$
\ell / p = . \ga_{_{\ell,1}} \!  \ga_{_{\ell,2}} \cdots \ga_{_{\ell,s}} \cdots  ,
$$
$\ga_{_{\ell,j}} \in \{0,1\} $ with $\ga_{_{\ell,j}} =0$
infinitely often. Hence, $\ell/p = \sum_{j=1}^\infty 2^{-j}
\ga_{_{\ell,j}}$. The fractional part of $\ell 2^{j} / p$ may be
written
$$
\{ \ell 2^{j} / p\} = .\ga_{_{\ell, j+1}}\ga_{_{\ell, j+2}} \cdots  \ga_{_{\ell,s}} \cdots  .
$$
Notice that  $\cos \frac{2\pi \ell 2^{j}}{ p} \leq 0$ if the
 fractional part of $\ell 2^{j} / p$ is (inclusively) between
 $1/4$ and $3/4$, which  follows if $\ga_{_{j+1}} \not=
 \ga_{_{j+2}}$. Thus, $\phi_s (\ell) $ is at least as large as the number of
 alterations in the sequence  $(\ga_{_{\ell,1}}, \ga_{_{\ell,2}},
 ..., \ga_{_{\ell, s+1}})$.

We now take $m$ such that $2^{m-1} < p < 2^m$. Observe that, for
$\ell=1, ..., p-1$, the subsequences $\ga (\ell) :=
(\ga_{_{\ell,1}}, \ga_{_{\ell,2}}, ..., \ga_{_{\ell, m}})$ of
length $m$ are pairwise distinct: If $\ga(\ell)=\ga(\ell')$ for
some $\ell< \ell'$ then $\frac{ \ell'-\ell}{p}$ is less than
$\sum_{j\geq m+1} 2^{-j} \leq 2^{-m}$, which is impossible as $p <
2^m$. Similarly, for fixed $j$ and $\ell=1,..., p-1$, all 
subsequences $\ga(\ell;j ) := 
(\ga_{_{\ell,j+1}}, \ga_{_{\ell,j+2}}, ..., \ga_{_{\ell, j+m}})$
are pairwise distinct. In
particular, for fixed $r$ with $r=0, ..., \lfloor s/m \rfloor -1$,
all subsequences $\ga(\ell;rm) $, $\ell=1,..., p-1$, are pairwise
distinct. Since the fractional part $\{ \frac{2^{rm} \ell }{p} \}=
.\ga_{_{\ell,rm+1}} \ga_{_{\ell,rm+2}} \cdots$ must be the same as
$\frac{ \ell' }{p}$ for some $\ell'$ in the range $1\leq \ell'
\leq p-1$, there is a unique permutation $\gs_r$ of $1, ... p-1$
such that $\ga(\ell;rm)= \ga ( \gs_r (\ell))$. Writing  $|\ga (
\gs_r (\ell))|_{_A} $ for the number of alternations in $\ga (
\gs_r (\ell))$, we have
$$ 
\phi_s (\ell ) \geq \sum_{r=0}^{\lfloor s/m \rfloor -1} |\ga ( \gs_r (\ell))|_{_A} , 
$$ 
where $\gs_{_0}$ is the identity. Therefore, \raf{phibd} gives
$$
\sum_{\ell=1}^{p-1} \Big|\prod_{j=0}^{s-1} \hat{\mu_j}
(\ell) \Big|^2 \leq \sum_{\ell=1}^{p-1}
 \, \xi^{2 \sum_{r=0}^{\lfloor s/m \rfloor -1} |\ga ( \gs_r
(\ell))|_{_A}}. 
$$

Using
\begin{eqnarray*}
\lefteqn{\xi^{x+y}+  \xi^{x'+y'}} \\
 &\leq& \xi^{\min\{x,x'\} +\min\{y,y'\}} +
 \xi^{\max\{x,x'\} +\max\{y,y'\}}
\end{eqnarray*}
inductively, the above upper bound may be maximized when all $\gs_r$'s  are  the identity, i.e.,
$$
\sum_{\ell=1}^{p-1} \Big|\prod_{j=0}^{s-1} \hat{\mu_j}
(\ell) \Big|^2 \leq \sum_{\ell=1}^{p-1}
\, \xi^{2\lfloor s/m \rfloor |\ga ( \ell)|_{_A}}.
$$
Note that $1/p \leq \ell /p \leq 1- 1/p$ implies that $\ga(\ell)$
is neither $(0,..., 0)$ nor $(1, ..., 1)$ (both are of length
$m$). This means that all $\ga(\ell)$ have at least one
alternation. Since $\ga(\ell)$'s are pairwise distinct,
$$
\sum_{\ell=1}^{p-1} \, \xi^{2\lfloor s/m \rfloor |\ga (
\ell)|_{_A}} \leq \sum_{\ga: |\ga|_{_A} >0} \xi^{2 \lfloor s/m
\rfloor|\ga|_{_A}},
$$
where the sum is taken over all sequences $\ga \in \{ 0, 1\}^m$
with $ |\ga|_{_A} >0$.

Let $H(z) $ be the number of $\ga$'s with exactly $z$ alterations.
Then
$$ H(z) =2 \binom{m-1}{z} ,$$
and hence
\begin{eqnarray*}
\sum_{\ga: |\ga|_{_A} >0} \xi^{2 \lfloor s/m \rfloor |\ga|_{_A}} 
  &=& 2 \sum_{z=1}^{m-1} \binom{m-1}{z} \xi^{2 \lfloor s/m \rfloor z} \\
  &=& 2 \left((1+\xi^{2\lfloor s/m \rfloor })^{m-1} -1 \right).
\end{eqnarray*}
\end{proof}

\begin{remark} \label{rmk:L2-to-sep}
To show a bound on the separation distance, we use Cauchy-Schwartz:
\begin{eqnarray*}
\lefteqn{\left| \frac{\P^{2s}(x,y)-\pi(y)}{\pi(y)}\right|^2} \\
    &=& \left| \frac{\sum_z \left(\P^s(x,z)-\pi(z)\right)\left(\P^s(z,y)-\pi(y)\right)}{\pi(y)} \right|^2 \\
    &=& \left| \sum_z \pi(z)\,\left(\frac{\P^s(x,z)}{\pi(z)}-1\right)\left(\frac{\P^{*s}(y,z)}{\pi(z)}-1\right)\right|^2 \\
\lefteqn{\leq \sum_z \pi(z)\,\left|\frac{\P^s(x,z)}{\pi(z)}-1\right|^2\,\sum_w \pi(w)\left|\frac{\P^{*s}(y,w)}{\pi(w)}-1\right|^2}
\end{eqnarray*}
For the ``block walk'' the first sum after the inequality is equal to the quantity upper bounded in equation \eqref{eqn:L2}, while the second is the same quantity but for the time-reversed walk $\P^*(a,b)=\pi(b)\P(b,a)/\pi(a)$. To bound the mixing time of the reversed walk let $b_i^*$ denote the sum of steps taken by $\R^*$ between the $(i-1)$-st and $i$th time that $j\to j/2$ is chosen (i.e. step size taken by time-reversed block walk), let $X_s^*=2^{-s+1}\,b_1^*+\cdots+b_s^*$ and let $b_i=-b_i^*$. Then 
\begin{eqnarray*}
\Pr[-2^{s-1}X_s^*=j] &=& \Pr[b_1+2b_2+\cdots+2^{s-1}b_s=j] \\
  &=& \Pr[X_s=j]
\end{eqnarray*}
because the $b_i$ are independent random variables from the same distribution as the blocks of $\R$. It follows from \eqref{eqn:L2} that
$$
1-\frac{\Pr[X_{2s}=y]}{u(y)}
   \leq 2 \left( (1+\xi^{2\lfloor s/m \rfloor})^{m-1} -1 \right)\,,
$$
and so after $2s\approx m\,\log_{\xi}\frac{\epsilon/2}{m-1}\leq 2m\log\frac{2(m-1)}{\epsilon}$ blocks the separation distance drops to $\epsilon$. 
\end{remark}

\begin{remark}
In recent work with Yuval Peres to appear in the journal version of this paper, we build on techniques from this paper and an idea from \cite{LPS03} and manage to improve the collision time bound to the conjectured value of $\Theta(\sqrt{p})$. The argument is based on showing that in $\Theta(\sqrt{p})$ steps the number of collisions $S$ in the $X_s$ walk satisfies $E(S)^2=\Theta(E(S^2))$, and so with constant probability there is a collision. This was in turn shown by re-writing $E(S^2)$ in terms of a quantity closely related to the Plancherel identity appearing in our Fourier proof of mixing time.
\end{remark}


\section{Concluding Remarks}

We sketch here some reasoning for why many common methods for bounding mixing times will not be useful in showing the optimal mixing bound in separation distance for the Pollard Rho walk.

A coupling argument bounds only the weaker total-variation distance, i.e. shows only that 
$$
\max_{A\subset G} \pi(A)-\Pr[X_t\in A] \leq \epsilon\,.
$$
To bound $\tau_s(1/2)$ with this requires $\epsilon\leq 1/2p$, which typically increases the mixing bound by a multiplicative factor of $\log(1/\epsilon)$. Total variation mixing time $\tau(1/2)$ is trivially at least $\log_3 p - 1$, and so this gives a separation bound of at best $O(\log^2 p)$. Alternatively, re-working the collision time argument in terms of variation distance results in a $\sqrt{\log p}$ loss.

When working with spectral gap, spectral profile, log-Sobolev and Nash inequalities a weakness arises in that mixing bounds in terms of these quantities are based on studying the rate of decay of variance. As such these do a poor job of distinguishing mixing time of a non-reversible walk $\P$ from that of its additive reversibilization $\P'=\frac{\P+\P^*}{2}$, or lazy additive reversibilization $\P''=\frac{\I}{2}+\frac{\P+\P^*}{4}$. The lazy additive reversibilization $\R''=\frac{\I}{2}+\frac{\R+\R^*}{4}$ of the Pollard Rho walk mixes in time $\tau_{s,R''}(1/2)=\Omega(\log^2 p)$ (see below), and so we expect that the aforementioned methods for bounding mixing time of $\R$ will do no better than this.

More precisely, the mixing time bounds involving these quantities can be shown by using the relation $\var(k_x^{t+1})-\var(k_x^t)=-\EE_{\P\P^*}(k_x^t,k_x^t)$ for the $t$-step density $k_x^t$ of a walk started at state $x$.
By a comparison argument it can be shown that if $k=p-1$ and $m\geq 2$ then $\frac{64m^2}{3}\,\EE_{(\R'')(\R'')^*}(f,f) \geq \EE_{\R^m\R^{*m}}(f,f)$. Hence these Dirichlet form based methods will show a mixing time bound on $\R^m$ which is no better than $\frac{64m^2}{3}$ times faster than the corresponding upper bound on the mixing time of $\R''$. To bound the mixing time of $\R''$ let $n=\frac 12\left\lfloor\log_2 p\right\rfloor$, $x=1/2^n\in \ZZ_p$ and $T=\frac{1}{9*32}\left(\left\lfloor \log_2 p\right\rfloor\right)^2$. It can then be shown that $(\R'')^T(x,S)\geq 7/9$ where
$$
S=\bigcup_{i=-3\sqrt{T/2}}^{3\sqrt{T/2}} \ \ 
\bigcup_{j=-T\,2^{n+3\sqrt{T/2}-1}}^{T\,2^{n+3\sqrt{T/2}-1}}
\left\{\frac{1+j}{2^{n-i}}\right\}\,.
$$
But $\pi(S)\leq \frac 1p\,(6\sqrt{T/2}+1)(T\,2^{n+3\sqrt{T/2}}+1)\leq 1/8$ and so for some $y\notin S$ we must have
$$
\frac{(\R'')^T(x,y)}{\pi(y)} \leq \frac{(\R'')^T(x,S^c)}{\pi(S^c)} 
< 3/10\,.
$$
It follows that $\tau_{s,\R'}(1/2)\geq T=\Omega(\log^2 p)$.  The corresponding mixing bound on $\R^m$ can then lead to a bound of at best
$$
\tau_{s,\R}(1/2)\leq m\tau_{s,\R^m}(1/2) = O((\log p)^2/m)\,.
$$
Hence if $m\ll \frac{\log p}{\log\log p}$ then none of these methods will match our $O(\log p\,\log \log p)$ mixing bound.

It thus appears that to show a better than $O(\log^2 p)$ mixing time bound it will be necessary to use a more specialized method, such as a more refined operator technique or computation involving a high power of the transition probability matrix. Two methods involving high powers were considered in this paper, a strong stationary time and a Fourier analysis argument. 

We now turn to limitations of the block approach of working with $X_s$ which was taken in our strong stationary time and Fourier analysis arguments. As with the random walk considered by \cite{AD86.1,CDG87.1}, we might expect
that the correct order of the mixing time of the $X_s$ walk considered in this paper is indeed $\Theta(\log p \log\log p)$, at least for $p$ of the form $2^m - 1$ and $k=p-1$. This is in fact the case by an argument fairly similar to that of Section 4 ``A proof of Case 2'' of \cite{Hil06.1}, which in turn closely follows a proof of \cite{CDG87.1}. The basic idea is by now fairly standard: choose a function and show that its expectation under the stationary distribution and under the $n$-step distribution $P^n$ are far apart, with sufficiently small variance to conclude that the two distributions ($P^n$ and $\pi$) must differ significantly.

In keeping with notation of \cite{Hil06.1}, suppose $p=2^t-1$ and let $k$ denote a variable over $\ZZ$ (no longer the exponent $y=x^k$). The ``separating function'' of interest $f:\ZZ_p\to\CC$ in this case is
$$
f(k) := \sum_{j=0}^{t-1} q^{k2^j}\qquad\textrm{where}\quad q=e^{2\pi i/p}\,.
$$
Then $E_U(f)=0$ if $p>1$, and $E_U(f\bar{f})=t$ and so $\var(f)=t$ where $\pi=U$ denotes the uniform distribution.
As in \cite{Hil06.1}, if $n=rt$ and $P_n$ denotes the distribution after $n$ steps of the block walk (i.e. $Y_0+X_n$) (the analysis uses $r=\delta\log t-d\in\NN$ for some fixed $\delta$), then $E_{P_n}(f) = t\,\Pi_1^r$ and $E_{P_n}(f\,\bar{f})  =   t\,\sum_{j=0}^{t-1}\Pi_j^r$ where $\Pi_j = \hat{P_t}(2^j-1)$, and so ${\rm Var}_{P_n}(f)=t\,\sum_{j=0}^{t-1} \Pi_j^r - t^2|\Pi_1|^{2r}$.

In order to bound variance and expectation we must approximate the $\Pi_j$. To do this, recall that $X_i=2X_{i-1}+b_i$; a generic such increment will be denoted by $b$, since the $b_i$ are independent random variables from the same distribution. Let $a_k=\Pr[b=k]=\Pr[b=-k]$. This satisfies the recurrence relation
$$
a_k = \frac 13\,(a_{k-1}+a_{k+1}),\,\quad a_0=\frac 13+\frac 23\,a_1,\quad a_{\infty}=0
$$
which can be solved to find that $a_k=\frac{1}{\sqrt{5}}\left(\frac{3-\sqrt{5}}{2}\right)^{|k|}$. It will also be useful to introduce a bit of notation. If $0\leq j\leq t-1$ then define $\mu_{\alpha}(x)=\Pr[b=x2^{-\alpha}]=a_{x 2^{-\alpha}}$ and
$$
\Pi_j = \hat{P_t}(2^j-1) = \prod_{\alpha=0}^{t-1} G\left(\frac{2^{\alpha}(2^j-1)}{p}\right)
$$
where
$$
G(x)=\frac{1}{\sqrt{5}}\,
     \frac{1-\left(\frac{3-\sqrt{5}}{2}\right)^2}
          {1+\left(\frac{3-\sqrt{5}}{2}\right)^2-(3-\sqrt{5})\cos(2\pi\,x)}\,.
$$

The remainder of the argument differs little from that of \cite{Hil06.1}. There is a small mistake in the proof of Claim 1 in \cite{Hil06.1}, but it does not effect the proof for the Rho walk.


\end{document}